\documentclass{amsart}
\usepackage{eurosym}
\usepackage{amsfonts}

\numberwithin{equation}{section}
\newtheorem{thm}{Theorem}[section]
\newtheorem{defn}{Definition}[section]

\newtheorem{lem}{Lemma}[section]
\newtheorem{crly}{Corollary}[section]

\begin{document}
\title[FRACTIONAL INTEGRAL INEQUALITIES FOR CONVEX AND CONCAVE FUNCTIONS]{%
FRACTIONAL INTEGRAL INEQUALITIES VIA ATANGANA-BALEANU OPERATORS FOR CONVEX
AND CONCAVE FUNCTIONS}
\author[A.O. Akdemir]{Ahmet Ocak Akdemir$^1$}
\address{$^1 $DEPARTMENT OF MATHEMATICS, FACULTY OF SCIENCES AND ARTS, A\v{G	%
}RI UNIVERSITY, A\v{G}RI, TURKEY}
\email{aocakakdemir@gmail.com}
\author[A. Karao\v{g}lan]{Ali Karao\v{g}lan$^2$}
\address{$^2$DEPARTMENT OF MATHEMATICS, FACULTY OF SCIENCES AND ARTS, ORDU
UNIVERSITY, ORDU, TURKEY}
\email{alikaraoglan@odu.edu.tr}
\author[M.A. Ragusa]{Maria Alessandra Ragusa$^3$}
\address{$^3$DIPARTIMENTO DI MATEMATICA E INFORMATICA, UNIVERSIT\'A DI CATANIA, CATANIA, ITALY - RUDN UNIVERSITY, MOSCOW, RUSSIA}
\email{mariaalessandra.ragusa@unict.it}
\author[E. Set]{Erhan Set$^{4}$ }
\address{$^1$DEPARTMENT OF MATHEMATICS, FACULTY OF SCIENCES AND ARTS, ORDU
UNIVERSITY, ORDU, TURKEY}
\email{erhanset@yahoo.com}
\dedicatory{}
\maketitle

\begin{abstract}
Recently, many fractional integral operators were introduced by different mathematicians. One of these fractional operators, Atangana-Baleanu fractional integral operator, was defined by Atangana and Baleanu in \cite{ABC}. In this study, firstly, a new identity by using Atangana-Baleanu fractional integral operators are proved. Then, new fractional integral inequalities have been obtained for convex and concave functions with the help of this identity and some certain integral inequalities.
\end{abstract}





\section{Introduction}

Mathematics is a tool that serves pure and applied sciences with its
deep-rooted history as old as human history, and sheds light on how to
express and then solve problems. Mathematics uses various concepts and their
relations with each other while performing this task. By defining spaces and
algebraic structures built on spaces, mathematics creates structures that
contribute to human life and nature. The concept of function is one of the
basic structures of mathematics, and many researchers have focused on new
function classes and made efforts to classify the space of functions. One of
the types of functions defined as a product of this intense effort is the
convex function, which has applications in statistics, inequality theory,
convex programming and numerical analysis. This interesting class of
functions is defined as follows.\newline

\begin{defn}
The mapping $f:[\theta _{1},\theta _{2}]\subseteq \mathbb{R}\rightarrow
\mathbb{R}$, is said to be convex if
\begin{equation}
f(\lambda x+(1-\lambda )y)\leq \lambda f(x)+(1-\lambda )f(y)  \label{1}
\end{equation}%
is valid for all $x,y\in \lbrack \theta _{1},\theta _{2}]$ and $\lambda \in
\lbrack 0,1]$.
\end{defn}

Many inequalities have been obtained by using this unique function type and
varieties in inequality theory, which is one of the most used areas of
convex functions. We will continue by introducing the Hermite-Hadamard
inequality that generate limits on the mean value of a convex function and
the famous Bullen inequality as follows. \newline
Assume that $f:I\subseteq
\mathbb{R}
\rightarrow
\mathbb{R}
$ is a convex mapping defined on the interval $I$ of $\mathbb{R}$ where $%
\theta _{1}<b.$ The following statement;%
\begin{equation}
f\left( \frac{\theta _{1}+\theta _{2}}{2}\right) \leq \frac{1}{\kappa
_{2}-\theta _{1}}\int\limits_{\theta _{1}}^{\theta _{2}}f(x)dx\leq \frac{%
f(\theta _{1})+f(\theta _{2})}{2}  \label{HH}
\end{equation}%
holds and known as Hermite-Hadamard inequality. Both inequalities hold in
the reversed direction if $f$ is concave. \newline
The Bullen's integral inequality can be presented as
\begin{equation*}
\frac{1}{\theta _{2}-\theta _{1}}\int_{\theta _{1}}^{\theta _{2}}f\left(
x\right) dx\leq \frac{1}{2}\left[ f\left( \frac{\theta _{1}+\theta _{2}}{2}%
\right) +\frac{f\left( \theta _{1}\right) +f\left( \theta _{2}\right) }{2}%
\right]
\end{equation*}%
where $f:I\subset
\mathbb{R}
\rightarrow
\mathbb{R}
$ is a convex mapping on the interval $I$ of $\mathbb{R}$ where $\kappa
_{1},\theta _{2}\in I$ with $\theta _{1}<\theta _{2}$.\newline
To provide detail information on convexity, let us consider some earlier
studies that have been performed by many researchers. In \cite{c6}, Jensen
introduced the concept of convex function to the literature for the first
time and drew attention to the fact that it seems to be the basis of the
concept of incremental function. In \cite{c8}, Beckenbach has mentioned
about the concept of convexity and emphasized several features of this
useful function class. In \cite{c7}, the authors have focused the relations
between convexity and Hermite-Hadamard's inequality. This study has led many
researchers to the link between convexity and integral inequalities, which
has guided studies in this field. Based on these studies, many papers have
been produced for different kinds of convex functions. In \cite{c2}, Akdemir
et al. have proved several new integral inequalities for
geometric-arithmetic convex functions via a new integral identity. Several
new Hadamard's type integral inequalities have been established with
applications to special means by Kavurmaci et al. in \cite{c3}. Therefore, a
similar argument has been carried out by Zhang et al. but now for $s-$%
geometrically convex functions in \cite{c5}. On all of these, Xi et al. have
extended the challenge to $m-$ and $(\alpha ,m)-$convex functions by
providing Hadamard type inequalities in \cite{c4}.\newline
Although fractional analysis has been known since ancient times, it has
recently become a more popular subject in mathematical analysis and applied
mathematics. The adventure that started with the question of whether the
solution will exist if the order is fractional in a differential equation,
has developed with many derivative and integral operators. By defining the
derivative and integral operators in fractional order, the researchers who
aimed to propose more effective solutions to the solution of physical
phenomena have turned to new operators with general and strong kernels over
time. This orientation has provided mathematics and applied sciences several
operators with kernel structures that differ in terms of locality and
singularity, as well as generalized operators with memory effect properties.
The struggle that started with the question of how the order in the
differential equation being a fraction would have consequences, has now
evolved into the problem of how to explain physical phenomena and find the
most effective fractional operators that will provide effective solutions to
real world problems. Let us introduce some fractional derivative and
integral operators that have broken ground in fractional analysis and have
proven their effectiveness in different fields by using by many researchers.%
\newline
We will remember the Caputo-Fabrizio derivative operators. Also, we would like to note that the functions belong to Hilbert spaces denoted by $H^{1}(0,\theta _{2})$

\begin{defn}
\cite{CF} Let $f\in H^{1}(0,\theta _{2})$, $\theta _{2}>\theta _{1}$, $%
\alpha \in \lbrack 0,1]$ then, the definition of the new Caputo fractional
derivative is:
\begin{equation}
^{CF}D^{\alpha }f(\tau_{1} )=\frac{M(\alpha )}{1-\alpha }\int_{\kappa
_{1}}^{\tau_{1} }f^{\prime }(s)exp\left[ -\frac{\alpha }{(1-\alpha )}%
(\tau_{1} -s)\right] ds  \label{tymfg89}
\end{equation}%
where $M(\alpha )$ is normalization function.
\end{defn}

Depends on this interesting fractional derivative operator, the authors have
defined the Caputo-Fabrizio fractional integral operator as follows.

\begin{defn}
\cite{AD} Let $f\in H^{1}(0,\theta _{2})$, $\theta _{2}>\theta _{1}$, $%
\alpha \in \lbrack 0,1]$ then, the definition of the left and right side of
Caputo-Fabrizio fractional integral is:
\begin{equation*}
\left( _{:::\theta _{1}}^{CF}I^{\alpha }\right) (\tau_{1} )=\frac{1-\alpha }{%
B(\alpha )}f(\tau_{1} )+\frac{\alpha }{B(\alpha )}\int_{\theta
_{1}}^{\tau_{1} }f(y)dy,
\end{equation*}%
and
\begin{equation*}
\left( ^{CF}I_{\theta _{2}}^{\alpha }\right) (\tau_{1} )=\frac{1-\alpha }{%
B(\alpha )}f(\tau_{1} )+\frac{\alpha }{B(\alpha )}\int_{\tau_{1} }^{\kappa
_{2}}f(y)dy
\end{equation*}%
where $B(\alpha )$ is normalization function.
\end{defn}

The Caputo-Fabrizio fractional derivative, which is used in dynamical systems,
physical phenomena, disease models and many other fields, is a highly
functional operator by definition, but has a deficiency in terms of not
meeting the initial conditions in the special case $\alpha=1$. The
improvement to eliminate this deficiency has been provided by the new
derivative operator developed by Atangana-Baleanu, which has versions in the
sense of Caputo and Riemann. In the sequel of this paper, we will denote the
normalization function with $B(\alpha)$ that the same properties with the $%
M(\alpha)$ which defined in Caputo-Fabrizio definition.

\begin{defn}
\cite{ABC} Let $f\in H^{1}(\theta _{1},\theta _{2})$, $\theta _{2}>\kappa
_{1}$, $\alpha \in \lbrack 0,1]$ then, the definition of the new fractional
derivative is given:
\begin{equation}
_{::::::\theta _{1}}^{ABC}D_{\tau_{1} }^{\alpha }\left[ f(\tau_{1} )\right] =%
\frac{B(\alpha )}{1-\alpha }\int_{a}^{\tau_{1} }f^{\prime }(x)E_{\alpha }%
\left[ -\alpha \frac{(\tau_{1} -x)^{\alpha }}{(1-\alpha )}\right] dx.
\label{abc}
\end{equation}
\end{defn}

\begin{defn}
\cite{ABC} Let $f\in H^{1}(\theta _{1},\theta _{2})$, $\theta _{2}>\kappa
_{1}$, $\alpha \in \lbrack 0,1]$ then, the definition of the new fractional
derivative is given:
\begin{equation}
_{:::::::\theta _{1}}^{ABR}D_{\tau_{1} }^{\alpha }\left[ f(\tau_{1} )\right]
=\frac{B(\alpha )}{1-\alpha }\frac{d}{d\tau_{1} }\int_{\theta
_{1}}^{\tau_{1} }f(x)E_{\alpha }\left[ -\alpha \frac{(\tau_{1} -x)^{\alpha }%
}{(1-\alpha )}\right] dx.  \label{abr}
\end{equation}%
\qquad
\end{defn}

Equations (\ref{abc}) and (\ref{abr}) have a non-local kernel. Also in
equation (\ref{abr}) when the function is constant we get zero.

The associated fractional integral operator has been defined by
Atangana-Baleanu as follows.

\begin{defn}
\cite{ABC} The fractional integral associate to the new fractional
derivative with non-local kernel of a function $f\in H^{1}(\kappa
_{1},\theta _{2})$ as defined:
\begin{equation*}
_{::::\theta _{1}}^{AB}I_{{}}^{\alpha }\left\{ f(\tau_{1} )\right\} =\frac{%
1-\alpha }{B(\alpha )}f(\tau_{1} )+\frac{\alpha }{B(\alpha )\Gamma (\alpha )}%
\int_{\theta _{1}}^{\tau_{1} }f(y)(\tau_{1} -y)^{\alpha -1}dy
\end{equation*}%
where $\theta _{2}>\theta _{1},\alpha \in \lbrack 0,1].$
\end{defn}

In \cite{TD}, Abdeljawad and Baleanu introduced right hand side of integral
operator as following; the right fractional new integral with ML kernel of
order $\alpha \in \lbrack 0,1]$ is defined by
\begin{equation*}
\left( ^{AB}I_{\theta _{2}}^{\alpha }\right) \left\{ f(\tau_{1} )\right\} =%
\frac{1-\alpha }{B(\alpha )}f(\tau_{1} )+\frac{\alpha }{B(\alpha )\Gamma
(\alpha )}\int_{\tau_{1} }^{\theta _{2}}f(y)(y-\tau_{1} )^{\alpha -1}dy.
\end{equation*}%
In \cite{AD}, Abdeljawad and Baleanu has presented some new results based on
fractional order derivatives and their discrete versions. Conformable
integral operators have been defined by Abdeljawad in \cite{13}. This useful
operator has been used to prove some new integral inequalities in \cite{15}.
Another important fractional operator -Riemann-Liouville fractional integral
operators- have been used to provide some new Simpson type integral
inequalities in \cite{11}. Ekinci and \"{O}zdemir have proved several
generalizations by using Riemann-Liouville fractional integral operators in
\cite{a8} and the authors have established some similar results with this
operator in \cite{e1}. In \cite{b11}, Akdemir et al. have presented some new
variants of celebrated Chebyshev inequality via generalized fractional
integral operators. The argument has been proceed with the study of Rashid
et al. (see \cite{a1}) that involves new investigations related to
generalized $k-$fractional integral operators. In \cite{a2}, Rashid et al. have presented some motivated findings that extend the argument to the Hilbert
spaces. For more information related to different kinds of
fractional operators, we recommend to consider \cite{8}. The applications of
fractional operators have been demonsrated by several researchers, we
suggest to see the papers \cite{e2}-\cite{e4}.\newline
The main motivation of this paper is to prove an integral identity that
includes the Atangana-Baleanu integral operator and to provide some new
Bullen type integral inequalities for differentiable convex and concave
functions with the help of this integral identity. Some special cases are
also considered.

\section{Main Results}

We will start with a new integral identity that will be used the proofs of
our main findings:

\begin{lem}
\label{ghd65} Let $f:[\theta _{1},\theta _{2}]\rightarrow \mathbb{R}$ be
differentiable function on $(\theta _{1},\theta _{2})$ with $\kappa
_{1}<\theta _{2}$. Then we have the following identity for Atangana-Baleanu
fractional integral operators
\begin{eqnarray*}
&&\frac{2(\theta _{2}-\theta _{1})^{\alpha }+(1-\alpha )2^{\alpha +1}\Gamma
(\alpha )}{(\theta _{2}-\theta _{1})^{\alpha +1}}\left[ f(\kappa
_{1})+f(\theta _{2})+2f\left( \frac{\theta _{1}+\theta _{2}}{2}\right) %
\right] \\
&&-\frac{2^{\alpha +1}B(\alpha )\Gamma (\alpha )}{(\theta _{2}-\kappa
_{1})^{\alpha +1}}\left[ ^{AB}I_{\frac{\theta _{1}+\theta _{2}}{2}}^{\alpha
}f(\theta _{1})+_{::::\theta _{1}}^{AB}I^{\alpha }f\left( \frac{\kappa
_{1}+\theta _{2}}{2}\right) \right. \\
&&\left. +_{::\frac{\theta _{1}+\theta _{2}}{2}}^{AB}I^{\alpha }f(\theta
_{2})+^{AB}I_{\theta _{2}}^{\alpha }f\left( \frac{\theta _{1}+\theta _{2}}{2}%
\right) \right] \\
&=&\int_{0}^{1}\left( (1-\tau _{1})^{\alpha }-\tau _{1}^{\alpha }\right)
f^{\prime }\left( \frac{1+\tau _{1}}{2}\theta _{1}+\frac{1-\tau _{1}}{2}%
\theta _{2}\right) d\tau _{1} \\
&&+\int_{0}^{1}\left( \tau _{1}^{\alpha }-(1-\tau _{1})^{\alpha }\right)
f^{\prime }\left( \frac{1+\tau _{1}}{2}\theta _{2}+\frac{1-\tau _{1}}{2}%
\kappa _{1}\right) d\tau _{1}
\end{eqnarray*}%
where $\alpha ,\tau _{1}\in \lbrack 0,1]$, $\Gamma (.)$ is Gamma function
and $B(\alpha )$ is normalization function.

\begin{proof}
By adding $I_{1}$ and $I_{2}$, we have
\begin{eqnarray*}
I_{1}+I_{2} &=&\int_{0}^{1}\left( (1-\tau _{1})^{\alpha }-\tau _{1}^{\alpha
}\right) f^{\prime }\left( \frac{1+\tau _{1}}{2}\theta _{1}+\frac{1-\tau _{1}%
}{2}\kappa _{2}\right) d\tau _{1} \\
&&+\int_{0}^{1}\left( \tau _{1}^{\alpha }-(1-\tau _{1})^{\alpha }\right)
f^{\prime }\left( \frac{1+\tau _{1}}{2}\theta _{2}+\frac{1-\tau _{1}}{2}%
\theta _{1}\right) d\tau _{1}.
\end{eqnarray*}%
By using integration, we have
\begin{eqnarray}
I_{1} &=&\int_{0}^{1}\left( (1-\tau _{1})^{\alpha }-\tau _{1}^{\alpha
}\right) f^{\prime }\left( \frac{1+\tau _{1}}{2}\theta _{1}+\frac{1-\tau _{1}%
}{2}\kappa _{2}\right) d\tau _{1}  \notag  \label{gh56+} \\
&=&\frac{\left( (1-\tau _{1})^{\alpha }-\tau _{1}^{\alpha }\right) f\left(
\frac{1+\tau _{1}}{2}\theta _{1}+\frac{1-\tau _{1}}{2}\theta _{2}\right)
d\tau _{1}}{\frac{\theta _{1}-\theta _{2}}{2}}\bigg|_{1}^{0}  \notag \\
&&-\frac{2\alpha }{\kappa _{2}-\theta _{1}}\int_{0}^{1}\left( (1-\tau
_{1})^{\alpha -1}+\tau _{1}^{\alpha -1}\right) f\left( \frac{1+\tau _{1}}{2}%
\theta _{1}+\frac{1-\tau _{1}}{2}\kappa _{2}\right) d\tau _{1}  \notag \\
&=&-\frac{2}{\theta _{1}-\theta _{2}}f(\theta _{1})-\frac{2}{\kappa
_{1}-\theta _{2}}f(\frac{\theta _{1}+\theta _{2}}{2})  \notag \\
&&-\frac{2\alpha }{\kappa _{2}-\theta _{1}}\int_{0}^{1}(1-\tau _{1})^{\alpha
-1}f\left( \frac{1+\tau _{1}}{2}\theta _{1}+\frac{1-\tau _{1}}{2}\theta
_{2}\right) d\tau _{1}  \notag \\
&&-\frac{2\alpha }{\theta _{2}-\theta _{1}}\int_{0}^{1}\tau _{1}^{\alpha
-1}f\left( \frac{1+\tau _{1}}{2}\theta _{1}+\frac{1-\tau _{1}}{2}\theta
_{2}\right) d\tau _{1}  \notag \\
&=&\frac{2}{\theta _{2}-\theta _{1}}\left( f(\theta _{1})+f(\frac{\kappa
_{1}+\theta _{2}}{2})\right) -\frac{2^{\alpha +1}\alpha }{(\kappa
_{2}-\theta _{1})^{\alpha +1}}\int_{\theta _{1}}^{\frac{\theta _{1}+\kappa
_{2}}{2}}\left( x-\theta _{1}\right) ^{\alpha -1}f(x)dx  \notag \\
&&-\frac{2^{\alpha +1}\alpha }{(\theta _{2}-\theta _{1})^{\alpha +1}}%
\int_{\theta _{1}}^{\frac{\theta _{1}+\theta _{2}}{2}}\left( \frac{\kappa
_{1}+\theta _{2}}{2}-x\right) ^{\alpha -1}f(x)dx.  \notag \\
&&
\end{eqnarray}%
Multiplying both side of \eqref{gh56+} identity by $\frac{(\kappa
_{2}-\theta _{1})^{\alpha +1}}{2^{\alpha +1}B(\alpha )\Gamma (\alpha )}$, we
have
\begin{eqnarray}
&&\frac{(\theta _{2}-\theta _{1})^{\alpha +1}}{2^{\alpha +1}B(\alpha )\Gamma
(\alpha )}I_{1}  \notag  \label{lk98n} \\
&=&\frac{(\theta _{2}-\theta _{1})^{\alpha }}{2^{\alpha }B(\alpha )\Gamma
(\alpha )}\left( f(\theta _{1})+f(\frac{\theta _{1}+\kappa _{2}}{2})\right) -%
\frac{\alpha }{B(\alpha )\Gamma (\alpha )}\int_{\kappa _{1}}^{\frac{\theta
_{1}+\theta _{2}}{2}}\left( x-\theta _{1}\right) ^{\alpha -1}f(x)dx  \notag
\\
&&-\frac{\alpha }{B(\alpha )\Gamma (\alpha )}\int_{\theta _{1}}^{\frac{%
\theta _{1}+\theta _{2}}{2}}\left( \frac{\theta _{1}+\theta _{2}}{2}%
-x\right) ^{\alpha -1}f(x)dx.  \notag \\
&&
\end{eqnarray}%
Similarly, by using integration, we get
\begin{eqnarray}
I_{2} &=&\int_{0}^{1}\left( \tau _{1}^{\alpha }-(1-\tau _{1})^{\alpha
}\right) f^{\prime }\left( \frac{1+\tau _{1}}{2}\theta _{2}+\frac{1-\tau _{1}%
}{2}\kappa _{1}\right) d\tau _{1}  \notag  \label{349f} \\
&=&\frac{\left( \tau _{1}^{\alpha }-(1-\tau _{1})^{\alpha }\right) f\left(
\frac{1+\tau _{1}}{2}\theta _{2}+\frac{1-\tau _{1}}{2}\theta _{1}\right)
d\tau _{1}}{\frac{\theta _{2}-\theta _{1}}{2}}\bigg|_{1}^{0}  \notag \\
&&-\frac{2\alpha }{\kappa _{2}-\theta _{1}}\int_{0}^{1}\left( \tau
_{1}^{\alpha -1}+(1-\tau _{1})^{\alpha -1}\right) f\left( \frac{1+\tau _{1}}{%
2}\theta _{2}+\frac{1-\tau _{1}}{2}\kappa _{1}\right) d\tau _{1}  \notag \\
&=&\frac{2}{\theta _{2}-\theta _{1}}\left( f(\theta _{2})+f(\frac{\kappa
_{1}+\theta _{2}}{2})\right) -\frac{2^{\alpha +1}\alpha }{(\kappa
_{2}-\theta _{1})^{\alpha +1}}\int_{\frac{\theta _{1}+\theta _{2}}{2}%
}^{\theta _{2}}\left( x-\frac{\theta _{1}+\theta _{2}}{2}\right) ^{\alpha
-1}f(x)dx  \notag \\
&&-\frac{2^{\alpha +1}\alpha }{(\theta _{2}-\theta _{1})^{\alpha +1}}\int_{%
\frac{\theta _{1}+\theta _{2}}{2}}^{\theta _{2}}\left( \theta _{2}-x\right)
^{\alpha -1}f(x)dx.  \notag \\
&&
\end{eqnarray}%
Multiplying both side of \eqref{349f} identity by $\frac{(\theta _{2}-\kappa
_{1})^{\alpha +1}}{2^{\alpha +1}B(\alpha )\Gamma (\alpha )}$, we get
\begin{eqnarray}
&&\frac{(\theta _{2}-\theta _{1})^{\alpha +1}}{2^{\alpha +1}B(\alpha )\Gamma
(\alpha )}I_{2}  \notag  \label{67dd0} \\
&=&\frac{(\theta _{2}-\theta _{1})^{\alpha }}{2^{\alpha }B(\alpha )\Gamma
(\alpha )}\left( f(\theta _{2})+f(\frac{\theta _{1}+\kappa _{2}}{2})\right) -%
\frac{\alpha }{B(\alpha )\Gamma (\alpha )}\int_{\frac{\theta _{1}+\theta _{2}%
}{2}}^{\theta _{2}}\left( x-\frac{\theta _{1}+\kappa _{2}}{2}\right)
^{\alpha -1}f(x)dx  \notag \\
&&-\frac{\alpha }{B(\alpha )\Gamma (\alpha )}\int_{\frac{\theta _{1}+\kappa
_{2}}{2}}^{\theta _{2}}\left( \theta _{2}-x\right) ^{\alpha -1}f(x)dx.
\notag \\
&&
\end{eqnarray}%
By adding identity \eqref{lk98n} and \eqref{67dd0}, we obtain
\begin{eqnarray*}
&&\frac{(\theta _{2}-\theta _{1})^{\alpha +1}}{2^{\alpha +1}B(\alpha )\Gamma
(\alpha )}[I_{1}+I_{2}] \\
&=&\frac{(\theta _{2}-\theta _{1})^{\alpha }+(1-\alpha )2^{\alpha }\Gamma
(\alpha )}{2^{\alpha }B(\alpha )\Gamma (\alpha )}\left[ f(\theta _{1})+f(%
\frac{\theta _{1}+\theta _{2}}{2})\right] \\
&&-\frac{1-\alpha }{B(\alpha )}f(\theta _{1})-\frac{\alpha }{B(\alpha
)\Gamma (\alpha )}\int_{\theta _{1}}^{\frac{\theta _{1}+\theta _{2}}{2}%
}\left( x-\theta _{1}\right) ^{\alpha -1}f(x)dx \\
&&-\frac{1-\alpha }{B(\alpha )}f\left( \frac{\theta _{1}+\theta _{2}}{2}%
\right) -\frac{\alpha }{B(\alpha )\Gamma (\alpha )}\int_{\theta _{1}}^{\frac{%
\theta _{1}+\theta _{2}}{2}}\left( \frac{\theta _{1}+\theta _{2}}{2}%
-x\right) ^{\alpha -1}f(x)dx \\
&+&\frac{(\theta _{2}-\theta _{1})^{\alpha }+(1-\alpha )2^{\alpha }\Gamma
(\alpha )}{2^{\alpha }B(\alpha )\Gamma (\alpha )}\left[ f(\theta _{2})+f(%
\frac{\theta _{1}+\theta _{2}}{2})\right] \\
&&-\frac{1-\alpha }{B(\alpha )}f(\theta _{2})-\frac{\alpha }{B(\alpha
)\Gamma (\alpha )}\int_{\frac{\kappa _{1}+\theta _{2}}{2}}^{\theta
_{2}}\left( \theta _{2}-x\right) ^{\alpha -1}f(x)dx \\
&&-\frac{1-\alpha }{B(\alpha )}f\left( \frac{\theta _{1}+\theta _{2}}{2}%
\right) -\frac{\alpha }{B(\alpha )\Gamma (\alpha )}\int_{\frac{\kappa
_{1}+\theta _{2}}{2}}^{\theta _{2}}\left( x-\frac{\theta _{1}+\theta _{2}}{2}%
\right) ^{\alpha -1}f(x)dx.
\end{eqnarray*}%
Using the definition of Atangana-Baleanu fractional integral operators, we
get
\begin{eqnarray*}
&&\frac{(\theta _{2}-\theta _{1})^{\alpha +1}}{2^{\alpha +1}B(\alpha )\Gamma
(\alpha )}\Bigg[\int_{0}^{1}\left( (1-\tau _{1})^{\alpha }-\tau _{1}^{\alpha
}\right) f^{\prime }\left( \frac{1+\tau _{1}}{2}\theta _{1}+\frac{1-\tau _{1}%
}{2}\theta _{2}\right) d\tau _{1} \\
&&+\int_{0}^{1}\left( \tau _{1}^{\alpha }-(1-\tau _{1})^{\alpha }\right)
f^{\prime }\left( \frac{1+\tau _{1}}{2}\theta _{2}+\frac{1-\tau _{1}}{2}%
\theta _{1}\right) d\tau _{1}\Bigg] \\
&=&\frac{(\theta _{2}-\theta _{1})^{\alpha }+(1-\alpha )2^{\alpha }\Gamma
(\alpha )}{2^{\alpha }B(\alpha )\Gamma (\alpha )}\left[ f(\kappa
_{1})+f(\theta _{2})+2f\left( \frac{\theta _{1}+\theta _{2}}{2}\right) %
\right] \\
&&-\left[ ^{AB}I_{\frac{\theta _{1}+\theta _{2}}{2}}^{\alpha }f(\kappa
_{1})+_{::::\theta _{1}}^{AB}I^{\alpha }f\left( \frac{\theta _{1}+\theta _{2}%
}{2}\right) +_{::\frac{\theta _{1}+\theta _{2}}{2}}^{AB}I^{\alpha }f(\theta
_{2})+^{AB}I_{\theta _{2}}^{\alpha }f\left( \frac{\theta _{1}+\theta _{2}}{2}%
\right) \right] .
\end{eqnarray*}
\end{proof}
\end{lem}

\begin{thm}
\label{ttr4h} Let $f:[\theta _{1},\theta _{2}]\rightarrow \mathbb{R}$ be
differentiable function on $(\theta _{1},\theta _{2})$ with $\kappa
_{1}<\theta _{2}$ and $f^{\prime }\in L_{1}[\theta _{1},\theta _{2}]$. If $%
|f^{\prime }|$ is a convex function, we have the following inequality for
Atangana-Baleanu fractional integral operators
\begin{eqnarray*}
&&\bigg|\frac{2(\theta _{2}-\theta _{1})^{\alpha }+(1-\alpha )2^{\alpha
+1}\Gamma (\alpha )}{(\theta _{2}-\theta _{1})^{\alpha +1}}\left[ f(\kappa
_{1})+f(\theta _{2})+2f\left( \frac{\theta _{1}+\theta _{2}}{2}\right) %
\right] \\
&&-\frac{2^{\alpha +1}B(\alpha )\Gamma (\alpha )}{(\theta _{2}-\kappa
_{1})^{\alpha +1}}\left[ ^{AB}I_{\frac{\theta _{1}+\theta _{2}}{2}}^{\alpha
}f(\theta _{1})+_{::::\theta _{1}}^{AB}I^{\alpha }f\left( \frac{\kappa
_{1}+\theta _{2}}{2}\right) \right. \\
&&\left. +_{::\frac{\theta _{1}+\theta _{2}}{2}}^{AB}I^{\alpha }f(\theta
_{2})+^{AB}I_{\theta _{2}}^{\alpha }f\left( \frac{\theta _{1}+\theta _{2}}{2}%
\right) \right] \bigg| \\
&\leq &\frac{2\left[ \left\vert f^{\prime }(\theta _{1})\right\vert
+\left\vert f^{\prime }(\theta _{2})\right\vert \right] }{\alpha +1}
\end{eqnarray*}%
where $\alpha \in \lbrack 0,1]$, $B(\alpha )$ is normalization function.

\begin{proof}
By using Lemma \ref{ghd65}, we can write
\begin{eqnarray*}
&&\bigg|\frac{2(\theta _{2}-\theta _{1})^{\alpha }+(1-\alpha )2^{\alpha
+1}\Gamma (\alpha )}{(\theta _{2}-\theta _{1})^{\alpha +1}}\left[ f(\kappa
_{1})+f(\theta _{2})+2f\left( \frac{\theta _{1}+\theta _{2}}{2}\right) %
\right] \\
&&-\frac{2^{\alpha +1}B(\alpha )\Gamma (\alpha )}{(\theta _{2}-\kappa
_{1})^{\alpha +1}}\left[ ^{AB}I_{\frac{\theta _{1}+\theta _{2}}{2}}^{\alpha
}f(\theta _{1})+_{::::\theta _{1}}^{AB}I^{\alpha }f\left( \frac{\kappa
_{1}+\theta _{2}}{2}\right) \right. \\
&&\left. +_{::\frac{\theta _{1}+\theta _{2}}{2}}^{AB}I^{\alpha }f(\theta
_{2})+^{AB}I_{\theta _{2}}^{\alpha }f\left( \frac{\theta _{1}+\theta _{2}}{2}%
\right) \right] \bigg| \\
&\leq &\int_{0}^{1}(1-\tau _{1})^{\alpha }\left\vert f^{\prime }\left( \frac{%
1+\tau _{1}}{2}\theta _{1}+\frac{1-\tau _{1}}{2}\theta _{2}\right)
\right\vert d\tau _{1} \\
&&+\int_{0}^{1}\tau _{1}^{\alpha }\left\vert f^{\prime }\left( \frac{1+\tau
_{1}}{2}\theta _{1}+\frac{1-\tau _{1}}{2}\theta _{2}\right) \right\vert
d\tau _{1} \\
&&+\int_{0}^{1}\tau _{1}^{\alpha }\left\vert f^{\prime }\left( \frac{1+\tau
_{1}}{2}\theta _{2}+\frac{1-\tau _{1}}{2}\theta _{1}\right) \right\vert
d\tau _{1} \\
&&+\int_{0}^{1}(1-\tau _{1})^{\alpha }\left\vert f^{\prime }\left( \frac{%
1+\tau _{1}}{2}\theta _{2}+\frac{1-\tau _{1}}{2}\theta _{1}\right)
\right\vert d\tau _{1}.
\end{eqnarray*}%
By using convexity of $|f^{\prime }|$, we get
\begin{eqnarray*}
&&\bigg|\frac{2(\theta _{2}-\theta _{1})^{\alpha }+(1-\alpha )2^{\alpha
+1}\Gamma (\alpha )}{(\theta _{2}-\theta _{1})^{\alpha +1}}\left[ f(\kappa
_{1})+f(\theta _{2})+2f\left( \frac{\theta _{1}+\theta _{2}}{2}\right) %
\right] \\
&&-\frac{2^{\alpha +1}B(\alpha )\Gamma (\alpha )}{(\theta _{2}-\kappa
_{1})^{\alpha +1}}\left[ ^{AB}I_{\frac{\theta _{1}+\theta _{2}}{2}}^{\alpha
}f(\theta _{1})+_{::::\theta _{1}}^{AB}I^{\alpha }f\left( \frac{\kappa
_{1}+\theta _{2}}{2}\right) \right. \\
&&\left. +_{::\frac{\theta _{1}+\theta _{2}}{2}}^{AB}I^{\alpha }f(\theta
_{2})+^{AB}I_{\theta _{2}}^{\alpha }f\left( \frac{\theta _{1}+\theta _{2}}{2}%
\right) \right] \bigg| \\
&\leq &\int_{0}^{1}(1-\tau _{1})^{\alpha }\left[ \frac{1+\tau _{1}}{2}%
\left\vert f^{\prime }(\theta _{1})\right\vert +\frac{1-\tau _{1}}{2}%
\left\vert f^{\prime }(\theta _{2})\right\vert \right] d\tau _{1} \\
&&+\int_{0}^{1}\tau _{1}^{\alpha }\left[ \frac{1+\tau _{1}}{2}\left\vert
f^{\prime }(\theta _{1})\right\vert +\frac{1-\tau _{1}}{2}\left\vert
f^{\prime }(\theta _{2})\right\vert \right] d\tau _{1} \\
&&+\int_{0}^{1}\tau _{1}^{\alpha }\left[ \frac{1+\tau _{1}}{2}\left\vert
f^{\prime }(\theta _{2})\right\vert +\frac{1-\tau _{1}}{2}\left\vert
f^{\prime }(\kappa _{1})\right\vert \right] d\tau _{1} \\
&&+\int_{0}^{1}(1-\tau _{1})^{\alpha }\left[ \frac{1+\tau _{1}}{2}\left\vert
f^{\prime }(\theta _{2})\right\vert +\frac{1-\tau _{1}}{2}\left\vert
f^{\prime }(\theta _{1})\right\vert \right] d\tau _{1}.
\end{eqnarray*}%
By computing the above integral, we obtain
\begin{eqnarray*}
&&\bigg|\frac{2(\theta _{2}-\theta _{1})^{\alpha }+(1-\alpha )2^{\alpha
+1}\Gamma (\alpha )}{(\theta _{2}-\theta _{1})^{\alpha +1}}\left[ f(\kappa
_{1})+f(\theta _{2})+2f\left( \frac{\theta _{1}+\theta _{2}}{2}\right) %
\right] \\
&&-\frac{2^{\alpha +1}B(\alpha )\Gamma (\alpha )}{(\theta _{2}-\kappa
_{1})^{\alpha +1}}\left[ ^{AB}I_{\frac{\theta _{1}+\theta _{2}}{2}}^{\alpha
}f(\theta _{1})+_{::::\theta _{1}}^{AB}I^{\alpha }f\left( \frac{\kappa
_{1}+\theta _{2}}{2}\right) \right. \\
&&\left. +_{::\frac{\theta _{1}+\theta _{2}}{2}}^{AB}I^{\alpha }f(\theta
_{2})+^{AB}I_{\theta _{2}}^{\alpha }f\left( \frac{\theta _{1}+\theta _{2}}{2}%
\right) \right] \bigg| \\
&\leq &\frac{2\left[ \left\vert f^{\prime }(\theta _{1})\right\vert
+\left\vert f^{\prime }(\theta _{2})\right\vert \right] }{\alpha +1}
\end{eqnarray*}%
and the proof is completed.
\end{proof}
\end{thm}

\begin{crly}
In Theorem \ref{ttr4h}, if we choose $\alpha =1$ we obtain
\begin{eqnarray*}
&&\bigg|\frac{f(\theta _{1})+f(\theta _{2})+2f\left( \frac{\theta
_{1}+\kappa _{2}}{2}\right) }{\theta _{2}-\theta _{1}}-\frac{4}{(\theta
_{2}-\kappa _{1})^{2}}\int_{\theta _{1}}^{\theta _{2}}f(x)dx\bigg| \\
&\leq &\frac{\left\vert f^{\prime }(\theta _{1})\right\vert +\left\vert
f^{\prime }(\kappa _{2})\right\vert }{2}.
\end{eqnarray*}
\end{crly}

\begin{thm}
\label{dfgre89} Let $f:[\theta _{1},\theta _{2}]\rightarrow \mathbb{R}$ be
differentiable function on $(\theta _{1},\theta _{2})$ with $\kappa
_{1}<\theta _{2}$ and $f^{\prime }\in L_{1}[\theta _{1},\theta _{2}]$. If $%
|f^{\prime }|^{q}$ is a convex function, then we have the following
inequality for Atangana-Baleanu fractional integral operators:
\begin{eqnarray*}
&&\bigg|\frac{2(\theta _{2}-\theta _{1})^{\alpha }+(1-\alpha )2^{\alpha
+1}\Gamma (\alpha )}{(\theta _{2}-\theta _{1})^{\alpha +1}}\left[ f(\kappa
_{1})+f(\theta _{2})+2f\left( \frac{\theta _{1}+\theta _{2}}{2}\right) %
\right] \\
&&-\frac{2^{\alpha +1}B(\alpha )\Gamma (\alpha )}{(\theta _{2}-\kappa
_{1})^{\alpha +1}}\left[ ^{AB}I_{\frac{\theta _{1}+\theta _{2}}{2}}^{\alpha
}f(\theta _{1})+_{::::\theta _{1}}^{AB}I^{\alpha }f\left( \frac{\kappa
_{1}+\theta _{2}}{2}\right) \right. \\
&&\left. +_{::\frac{\theta _{1}+\theta _{2}}{2}}^{AB}I^{\alpha }f(\theta
_{2})+^{AB}I_{\theta _{2}}^{\alpha }f\left( \frac{\theta _{1}+\theta _{2}}{2}%
\right) \right] \bigg| \\
&\leq &\frac{2}{(\alpha p+1)^{\frac{1}{p}}}\left[ \left( \frac{3\left\vert
f^{\prime }(\theta _{1})\right\vert ^{q}+\left\vert f^{\prime }(\kappa
_{2})\right\vert ^{q}}{4}\right) ^{\frac{1}{q}}+\left( \frac{3\left\vert
f^{\prime }(\theta _{2})\right\vert ^{q}+\left\vert f^{\prime }(\kappa
_{1})\right\vert ^{q}}{4}\right) ^{\frac{1}{q}}\right]
\end{eqnarray*}%
where $p^{-1}+q^{-1}=1$, $\alpha \in \lbrack 0,1]$, $q>1$, $B(\alpha )$ is
normalization function.

\begin{proof}
By using the identity that is given in Lemma \ref{ghd65}, we have
\begin{eqnarray*}
&&\bigg|\frac{2(\theta _{2}-\theta _{1})^{\alpha }+(1-\alpha )2^{\alpha
+1}\Gamma (\alpha )}{(\theta _{2}-\theta _{1})^{\alpha +1}}\left[ f(\kappa
_{1})+f(\theta _{2})+2f\left( \frac{\theta _{1}+\theta _{2}}{2}\right) %
\right] \\
&&-\frac{2^{\alpha +1}B(\alpha )\Gamma (\alpha )}{(\theta _{2}-\kappa
_{1})^{\alpha +1}}\left[ ^{AB}I_{\frac{\theta _{1}+\theta _{2}}{2}}^{\alpha
}f(\theta _{1})+_{::::\theta _{1}}^{AB}I^{\alpha }f\left( \frac{\kappa
_{1}+\theta _{2}}{2}\right) \right. \\
&&\left. +_{::\frac{\theta _{1}+\theta _{2}}{2}}^{AB}I^{\alpha }f(\theta
_{2})+^{AB}I_{\theta _{2}}^{\alpha }f\left( \frac{\theta _{1}+\theta _{2}}{2}%
\right) \right] \bigg| \\
&\leq &\int_{0}^{1}(1-\tau _{1})^{\alpha }\left\vert f^{\prime }\left( \frac{%
1+\tau _{1}}{2}\theta _{1}+\frac{1-\tau _{1}}{2}\theta _{2}\right)
\right\vert d\tau _{1}+\int_{0}^{1}\tau _{1}^{\alpha }\left\vert f^{\prime
}\left( \frac{1+\tau _{1}}{2}\theta _{1}+\frac{1-\tau _{1}}{2}\theta
_{2}\right) \right\vert d\tau _{1} \\
&&+\int_{0}^{1}\tau _{1}^{\alpha }\left\vert f^{\prime }\left( \frac{1+\tau
_{1}}{2}\theta _{2}+\frac{1-\tau _{1}}{2}\theta _{1}\right) \right\vert
d\tau _{1}+\int_{0}^{1}(1-\tau _{1})^{\alpha }\left\vert f^{\prime }\left(
\frac{1+\tau _{1}}{2}\theta _{2}+\frac{1-\tau _{1}}{2}\theta _{1}\right)
\right\vert d\tau _{1}.
\end{eqnarray*}%
By applying H\"{o}lder inequality, we have
\begin{eqnarray*}
&&\bigg|\frac{2(\theta _{2}-\theta _{1})^{\alpha }+(1-\alpha )2^{\alpha
+1}\Gamma (\alpha )}{(\theta _{2}-\theta _{1})^{\alpha +1}}\left[ f(\kappa
_{1})+f(\theta _{2})+2f\left( \frac{\theta _{1}+\theta _{2}}{2}\right) %
\right] \\
&&-\frac{2^{\alpha +1}B(\alpha )\Gamma (\alpha )}{(\theta _{2}-\kappa
_{1})^{\alpha +1}}\left[ ^{AB}I_{\frac{\theta _{1}+\theta _{2}}{2}}^{\alpha
}f(\theta _{1})+_{::::\theta _{1}}^{AB}I^{\alpha }f\left( \frac{\kappa
_{1}+\theta _{2}}{2}\right) \right. \\
&&\left. +_{::\frac{\theta _{1}+\theta _{2}}{2}}^{AB}I^{\alpha }f(\theta
_{2})+^{AB}I_{\theta _{2}}^{\alpha }f\left( \frac{\theta _{1}+\theta _{2}}{2}%
\right) \right] \bigg| \\
&\leq &\left( \int_{0}^{1}(1-\tau _{1})^{\alpha p}d\tau _{1}\right) ^{\frac{1%
}{p}}\left( \int_{0}^{1}\left\vert f^{\prime }\left( \frac{1+\tau _{1}}{2}%
\kappa _{1}+\frac{1-\tau _{1}}{2}\theta _{2}\right) \right\vert ^{q}d\tau
_{1}\right) ^{\frac{1}{q}} \\
&&+\left( \int_{0}^{1}\tau _{1}^{\alpha p}d\tau _{1}\right) ^{\frac{1}{p}%
}\left( \int_{0}^{1}\left\vert f^{\prime }\left( \frac{1+\tau _{1}}{2}\theta
_{1}+\frac{1-\tau _{1}}{2}\theta _{2}\right) \right\vert ^{q}d\tau
_{1}\right) ^{\frac{1}{q}} \\
&&+\left( \int_{0}^{1}\tau _{1}^{\alpha p}d\tau _{1}\right) ^{\frac{1}{p}%
}\left( \int_{0}^{1}\left\vert f^{\prime }\left( \frac{1+\tau _{1}}{2}\theta
_{2}+\frac{1-\tau _{1}}{2}\theta _{1}\right) \right\vert ^{q}d\tau
_{1}\right) ^{\frac{1}{q}} \\
&&+\left( \int_{0}^{1}(1-\tau _{1})^{\alpha p}d\tau _{1}\right) ^{\frac{1}{p}%
}\left( \int_{0}^{1}\left\vert f^{\prime }\left( \frac{1+\tau _{1}}{2}\theta
_{2}+\frac{1-\tau _{1}}{2}\theta _{1}\right) \right\vert ^{q}d\tau
_{1}\right) ^{\frac{1}{q}}.
\end{eqnarray*}%
By using convexity of $|f^{\prime }|^{q}$, we obtain
\begin{eqnarray*}
&&\bigg|\frac{2(\theta _{2}-\theta _{1})^{\alpha }+(1-\alpha )2^{\alpha
+1}\Gamma (\alpha )}{(\theta _{2}-\theta _{1})^{\alpha +1}}\left[ f(\kappa
_{1})+f(\theta _{2})+2f\left( \frac{\theta _{1}+\theta _{2}}{2}\right) %
\right] \\
&&-\frac{2^{\alpha +1}B(\alpha )\Gamma (\alpha )}{(\theta _{2}-\kappa
_{1})^{\alpha +1}}\left[ ^{AB}I_{\frac{\theta _{1}+\theta _{2}}{2}}^{\alpha
}f(\theta _{1})+_{::::\theta _{1}}^{AB}I^{\alpha }f\left( \frac{\kappa
_{1}+\theta _{2}}{2}\right) \right. \\
&&\left. +_{::\frac{\theta _{1}+\theta _{2}}{2}}^{AB}I^{\alpha }f(\theta
_{2})+^{AB}I_{\theta _{2}}^{\alpha }f\left( \frac{\theta _{1}+\theta _{2}}{2}%
\right) \right] \bigg| \\
&\leq &\left( \int_{0}^{1}(1-\tau _{1})^{\alpha p}d\tau _{1}\right) ^{\frac{1%
}{p}}\left( \int_{0}^{1}\left[ \frac{1+\tau _{1}}{2}\left\vert f^{\prime
}(\kappa _{1})\right\vert ^{q}+\frac{1-\tau _{1}}{2}\left\vert f^{\prime
}(\kappa _{2})\right\vert ^{q}\right] dv\right) ^{\frac{1}{q}} \\
&&+\left( \int_{0}^{1}\tau _{1}^{\alpha p}d\tau _{1}\right) ^{\frac{1}{p}%
}\left( \int_{0}^{1}\left[ \frac{1+\tau _{1}}{2}\left\vert f^{\prime
}(\kappa _{1})\right\vert ^{q}+\frac{1-\tau _{1}}{2}\left\vert f^{\prime
}(\kappa _{2})\right\vert ^{q}\right] d\tau _{1}\right) ^{\frac{1}{q}} \\
&&+\left( \int_{0}^{1}\tau _{1}^{\alpha p}d\tau _{1}\right) ^{\frac{1}{p}%
}\left( \int_{0}^{1}\left[ \frac{1+\tau _{1}}{2}\left\vert f^{\prime
}(\kappa _{2})\right\vert ^{q}+\frac{1-\tau _{1}}{2}\left\vert f^{\prime
}(\kappa _{1})\right\vert ^{q}\right] d\tau _{1}\right) ^{\frac{1}{q}} \\
&&+\left( \int_{0}^{1}(1-\tau _{1})^{\alpha p}d\tau _{1}\right) ^{\frac{1}{p}%
}\left( \int_{0}^{1}\left[ \frac{1+\tau _{1}}{2}\left\vert f^{\prime
}(\kappa _{2})\right\vert ^{q}+\frac{1-\tau _{1}}{2}\left\vert f^{\prime
}(\kappa _{1})\right\vert ^{q}\right] d\tau _{1}\right) ^{\frac{1}{q}}.
\end{eqnarray*}%
By calculating the integrals that is in the above inequalities, we get
desired result.
\end{proof}
\end{thm}

\begin{crly}
In Theorem \ref{dfgre89}, if we choose $\alpha =1$ we obtain
\begin{eqnarray*}
&&\bigg|\frac{f(\theta _{1})+f(\theta _{2})+2f\left( \frac{\kappa
_{1}+\theta _{2}}{2}\right) }{\theta _{2}-\theta _{1}}-\frac{4}{(\kappa
_{2}-\theta _{1})^{2}}\int_{\theta _{1}}^{\theta _{2}}f(x)dx\bigg| \\
&\leq &\frac{1}{(p+1)^{\frac{1}{p}}}\bigg[\left( \frac{3\left\vert f^{\prime
}(\theta _{1})\right\vert ^{q}+\left\vert f^{\prime }(\kappa
_{2})\right\vert ^{q}}{4}\right) ^{\frac{1}{q}} \\
&&+\left( \frac{3\left\vert f^{\prime }(\theta _{2})\right\vert
^{q}+\left\vert f^{\prime }(\theta _{1})\right\vert ^{q}}{4}\right) ^{\frac{1%
}{q}}\bigg].
\end{eqnarray*}
\end{crly}

\begin{thm}
\label{908tty6} Let $f:[\theta _{1},\theta _{2}]\rightarrow \mathbb{R}$ be
differentiable function on $(\theta _{1},\theta _{2})$ with $\kappa
_{1}<\theta _{2}$ and $f^{\prime }\in L_{1}[\theta _{1},\theta _{2}]$. If $%
|f^{\prime }|^{q}$ is a convex function, then we have the following
inequality for Atangana-Baleanu fractional integral operators
\begin{eqnarray*}
&&\bigg|\frac{2(\theta _{2}-\theta _{1})^{\alpha }+(1-\alpha )2^{\alpha
+1}\Gamma (\alpha )}{(\theta _{2}-\theta _{1})^{\alpha +1}}\left[ f(\kappa
_{1})+f(\theta _{2})+2f\left( \frac{\theta _{1}+\theta _{2}}{2}\right) %
\right] \\
&&-\frac{2^{\alpha +1}B(\alpha )\Gamma (\alpha )}{(\theta _{2}-\kappa
_{1})^{\alpha +1}}\left[ ^{AB}I_{\frac{\theta _{1}+\theta _{2}}{2}}^{\alpha
}f(\theta _{1})+_{::::\theta _{1}}^{AB}I^{\alpha }f\left( \frac{\kappa
_{1}+\theta _{2}}{2}\right) \right. \\
&&\left. +_{::\frac{\theta _{1}+\theta _{2}}{2}}^{AB}I^{\alpha }f(\theta
_{2})+^{AB}I_{\theta _{2}}^{\alpha }f\left( \frac{\theta _{1}+\theta _{2}}{2}%
\right) \right] \bigg| \\
&\leq &\left( \frac{1}{\alpha +1}\right) ^{1-\frac{1}{q}}\Bigg[\left( \frac{%
\alpha +3}{2(\alpha +1)(\alpha +2)}\left\vert f^{\prime }(\kappa
_{1})\right\vert ^{q}+\frac{1}{2(\alpha +2)}\left\vert f^{\prime }(\kappa
_{2})\right\vert ^{q}\right) ^{\frac{1}{q}} \\
&&+\left( \frac{2\alpha +3}{2(\alpha +1)(\alpha +2)}\left\vert f^{\prime
}(\theta _{1})\right\vert ^{q}+\frac{1}{2(\alpha +1)(\alpha +2)}\left\vert
f^{\prime }(\theta _{2})\right\vert ^{q}\right) ^{\frac{1}{q}} \\
&&+\left( \frac{2\alpha +3}{2(\alpha +1)(\alpha +2)}\left\vert f^{\prime
}(\theta _{2})\right\vert ^{q}+\frac{1}{2(\alpha +1)(\alpha +2)}\left\vert
f^{\prime }(\theta _{1})\right\vert ^{q}\right) ^{\frac{1}{q}} \\
&&+\left( \frac{\alpha +3}{2(\alpha +1)(\alpha +2)}\left\vert f^{\prime
}(\theta _{2})\right\vert ^{q}+\frac{1}{2(\alpha +2)}\left\vert f^{\prime
}(\theta _{1})\right\vert ^{q}\right) ^{\frac{1}{q}}\Bigg]
\end{eqnarray*}%
where $\alpha \in \lbrack 0,1]$, $q\geq 1$, $B(\alpha )$ is normalization
function.

\begin{proof}
By Lemma \ref{ghd65}, we get
\begin{eqnarray*}
&&\bigg|\frac{2(\theta _{2}-\theta _{1})^{\alpha }+(1-\alpha )2^{\alpha
+1}\Gamma (\alpha )}{(\theta _{2}-\theta _{1})^{\alpha +1}}\left[ f(\kappa
_{1})+f(\theta _{2})+2f\left( \frac{\theta _{1}+\theta _{2}}{2}\right) %
\right]  \\
&&-\frac{2^{\alpha +1}B(\alpha )\Gamma (\alpha )}{(\theta _{2}-\kappa
_{1})^{\alpha +1}}\left[ ^{AB}I_{\frac{\theta _{1}+\theta _{2}}{2}}^{\alpha
}f(\theta _{1})+_{::::\theta _{1}}^{AB}I^{\alpha }f\left( \frac{\kappa
_{1}+\theta _{2}}{2}\right) \right.  \\
&&\left. +_{::\frac{\theta _{1}+\theta _{2}}{2}}^{AB}I^{\alpha }f(\theta
_{2})+^{AB}I_{\theta _{2}}^{\alpha }f\left( \frac{\theta _{1}+\theta _{2}}{2}%
\right) \right] \bigg| \\
&\leq &\int_{0}^{1}(1-\tau _{1})^{\alpha }\left\vert f^{\prime }\left( \frac{%
1+\tau _{1}}{2}\theta _{1}+\frac{1-\tau _{1}}{2}\theta _{2}\right)
\right\vert d\tau _{1}+\int_{0}^{1}t^{\alpha }\left\vert f^{\prime }\left(
\frac{1+\tau _{1}}{2}\kappa _{1}+\frac{1-\tau _{1}}{2}\theta _{2}\right)
\right\vert d\tau _{1} \\
&&+\int_{0}^{1}\tau _{1}^{\alpha }\left\vert f^{\prime }\left( \frac{1+\tau
_{1}}{2}\theta _{2}+\frac{1-\tau _{1}}{2}\theta _{1}\right) \right\vert
d\tau _{1}+\int_{0}^{1}(1-\tau _{1})^{\alpha }\left\vert f^{\prime }\left(
\frac{1+\tau _{1}}{2}\theta _{2}+\frac{1-\tau _{1}}{2}\theta _{1}\right)
\right\vert d\tau _{1}.
\end{eqnarray*}%
By applying power mean inequality, we get
\begin{eqnarray*}
&&\bigg|\frac{2(\theta _{2}-\theta _{1})^{\alpha }+(1-\alpha )2^{\alpha
+1}\Gamma (\alpha )}{(\theta _{2}-\theta _{1})^{\alpha +1}}\left[ f(\kappa
_{1})+f(\theta _{2})+2f\left( \frac{\theta _{1}+\theta _{2}}{2}\right) %
\right]  \\
&&-\frac{2^{\alpha +1}B(\alpha )\Gamma (\alpha )}{(\theta _{2}-\kappa
_{1})^{\alpha +1}}\left[ ^{AB}I_{\frac{\theta _{1}+\theta _{2}}{2}}^{\alpha
}f(\theta _{1})+_{::::\theta _{1}}^{AB}I^{\alpha }f\left( \frac{\kappa
_{1}+\theta _{2}}{2}\right) \right.  \\
&&\left. +_{::\frac{\theta _{1}+\theta _{2}}{2}}^{AB}I^{\alpha }f(\theta
_{2})+^{AB}I_{\theta _{2}}^{\alpha }f\left( \frac{\theta _{1}+\theta _{2}}{2}%
\right) \right] \bigg| \\
&\leq &\left( \int_{0}^{1}(1-\tau _{1})^{\alpha }d\tau _{1}\right) ^{1-\frac{%
1}{q}}\left( \int_{0}^{1}(1-\tau _{1})^{\alpha }\left\vert f^{\prime }\left(
\frac{1+\tau _{1}}{2}\theta _{1}+\frac{1-\tau _{1}}{2}\theta _{2}\right)
\right\vert ^{q}d\tau _{1}\right) ^{\frac{1}{q}} \\
&&+\left( \int_{0}^{1}\tau _{1}^{\alpha }d\tau _{1}\right) ^{1-\frac{1}{q}%
}\left( \int_{0}^{1}\tau _{1}^{\alpha }\left\vert f^{\prime }\left( \frac{%
1+\tau _{1}}{2}\theta _{1}+\frac{1-\tau _{1}}{2}\theta _{2}\right)
\right\vert ^{q}d\tau _{1}\right) ^{\frac{1}{q}} \\
&&+\left( \int_{0}^{1}\tau _{1}^{\alpha }d\tau _{1}\right) ^{1-\frac{1}{q}%
}\left( \int_{0}^{1}\tau _{1}^{\alpha }\left\vert f^{\prime }\left( \frac{%
1+\tau _{1}}{2}\theta _{2}+\frac{1-\tau _{1}}{2}\theta _{1}\right)
\right\vert ^{q}d\tau _{1}\right) ^{\frac{1}{q}} \\
&&+\left( \int_{0}^{1}(1-\tau _{1})^{\alpha }d\tau _{1}\right) ^{1-\frac{1}{q%
}}\left( \int_{0}^{1}(1-\tau _{1})^{\alpha }\left\vert f^{\prime }\left(
\frac{1+\tau _{1}}{2}\theta _{2}+\frac{1-\tau _{1}}{2}\theta _{1}\right)
\right\vert ^{q}d\tau _{1}\right) ^{\frac{1}{q}}.
\end{eqnarray*}%
By using convexity of $|f^{\prime }|^{q}$, we obtain
\begin{eqnarray*}
&&\bigg|\frac{2(\theta _{2}-\theta _{1})^{\alpha }+(1-\alpha )2^{\alpha
+1}\Gamma (\alpha )}{(\theta _{2}-\theta _{1})^{\alpha +1}}\left[ f(\kappa
_{1})+f(\theta _{2})+2f\left( \frac{\theta _{1}+\theta _{2}}{2}\right) %
\right]  \\
&&-\frac{2^{\alpha +1}B(\alpha )\Gamma (\alpha )}{(\theta _{2}-\kappa
_{1})^{\alpha +1}}\left[ ^{AB}I_{\frac{\theta _{1}+\theta _{2}}{2}}^{\alpha
}f(\theta _{1})+_{::::\theta _{1}}^{AB}I^{\alpha }f\left( \frac{\kappa
_{1}+\theta _{2}}{2}\right) \right.  \\
&&\left. +_{::\frac{\theta _{1}+\theta _{2}}{2}}^{AB}I^{\alpha }f(\theta
_{2})+^{AB}I_{\theta _{2}}^{\alpha }f\left( \frac{\theta _{1}+\theta _{2}}{2}%
\right) \right] \bigg| \\
&\leq &\left( \int_{0}^{1}(1-\tau _{1})^{\alpha }d\tau _{1}\right) ^{1-\frac{%
1}{q}}\left( \int_{0}^{1}(1-\tau _{1})^{\alpha }\left[ \frac{1+\tau _{1}}{2}%
\left\vert f^{\prime }(\theta _{1})\right\vert ^{q}+\frac{1-\tau _{1}}{2}%
\left\vert f^{\prime }(\theta _{2})\right\vert ^{q}\right] d\tau _{1}\right)
^{\frac{1}{q}} \\
&&+\left( \int_{0}^{1}\tau _{1}^{\alpha }d\tau _{1}\right) ^{1-\frac{1}{q}%
}\left( \int_{0}^{1}\tau _{1}^{\alpha }\left[ \frac{1+\tau _{1}}{2}%
\left\vert f^{\prime }(\theta _{1})\right\vert ^{q}+\frac{1-\tau _{1}}{2}%
\left\vert f^{\prime }(\theta _{2})\right\vert ^{q}\right] d\tau _{1}\right)
^{\frac{1}{q}} \\
&&+\left( \int_{0}^{1}\tau _{1}^{\alpha }d\tau _{1}\right) ^{1-\frac{1}{q}%
}\left( \int_{0}^{1}\tau _{1}^{\alpha }\left[ \frac{1+\tau _{1}}{2}%
\left\vert f^{\prime }(\theta _{2})\right\vert ^{q}+\frac{1-\tau _{1}}{2}%
\left\vert f^{\prime }(\theta _{1})\right\vert ^{q}\right] d\tau _{1}\right)
^{\frac{1}{q}} \\
&&+\left( \int_{0}^{1}(1-\tau _{1})^{\alpha }d\tau _{1}\right) ^{1-\frac{1}{q%
}}\left( \int_{0}^{1}(1-\tau _{1})^{\alpha }\left[ \frac{1+\tau _{1}}{2}%
\left\vert f^{\prime }(\theta _{2})\right\vert ^{q}+\frac{1-\tau _{1}}{2}%
\left\vert f^{\prime }(\theta _{1})\right\vert ^{q}\right] d\tau _{1}\right)
^{\frac{1}{q}}
\end{eqnarray*}%
By computing the above integrals, the proof is completed.
\end{proof}
\end{thm}

\begin{crly}
In Theorem \ref{908tty6}, if we choose $\alpha =1$ we obtain
\begin{eqnarray*}
&&\bigg|\frac{2\left[ f(\theta _{1})+f(\theta _{2})+2f\left( \frac{\kappa
_{1}+\theta _{2}}{2}\right) \right] }{\theta _{2}-\theta _{1}}-\frac{8}{%
(\theta _{2}-\theta _{1})^{2}}\int_{\theta _{1}}^{\theta _{2}}f(x)dx\bigg| \\
&\leq &\left( \frac{1}{2}\right) ^{1-\frac{1}{q}}\Bigg[\left( \frac{%
2\left\vert f^{\prime }(\theta _{1})\right\vert ^{q}+\left\vert f^{\prime
}(\theta _{2})\right\vert ^{q}}{6}\right) ^{\frac{1}{q}}+\left( \frac{%
5\left\vert f^{\prime }(\theta _{1})\right\vert ^{q}+\left\vert f^{\prime
}(\theta _{2})\right\vert ^{q}}{12}\right) ^{\frac{1}{q}} \\
&&+\left( \frac{5\left\vert f^{\prime }(\theta _{2})\right\vert
^{q}+\left\vert f^{\prime }(\theta _{1})\right\vert ^{q}}{12}\right) ^{\frac{%
1}{q}}+\left( \frac{2\left\vert f^{\prime }(\theta _{2})\right\vert
^{q}+\left\vert f^{\prime }(\theta _{1})\right\vert ^{q}}{6}\right) ^{\frac{1%
}{q}}\Bigg].
\end{eqnarray*}
\end{crly}

\begin{thm}
\label{9888um} Let $f:[\theta _{1},\theta _{2}]\rightarrow \mathbb{R}$ be
differentiable function on $(\theta _{1},\theta _{2})$ with $\kappa
_{1}<\theta _{2}$ and $f^{\prime }\in L_{1}[\theta _{1},\theta _{2}].$ If $%
|f^{\prime }|^{q}$ is a convex function, then we have the following
inequality for Atangana-Baleanu fractional integral operators
\begin{eqnarray*}
&&\bigg|\frac{2(\theta _{2}-\theta _{1})^{\alpha }+(1-\alpha )2^{\alpha
+1}\Gamma (\alpha )}{(\theta _{2}-\theta _{1})^{\alpha +1}}\left[ f(\kappa
_{1})+f(\theta _{2})+2f\left( \frac{\theta _{1}+\theta _{2}}{2}\right) %
\right] \\
&&-\frac{2^{\alpha +1}B(\alpha )\Gamma (\alpha )}{(\theta _{2}-\kappa
_{1})^{\alpha +1}}\left[ ^{AB}I_{\frac{\theta _{1}+\theta _{2}}{2}}^{\alpha
}f(\theta _{1})+_{::::\theta _{1}}^{AB}I^{\alpha }f\left( \frac{\kappa
_{1}+\theta _{2}}{2}\right) \right. \\
&&\left. +_{::\frac{\theta _{1}+\theta _{2}}{2}}^{AB}I^{\alpha }f(\theta
_{2})+^{AB}I_{\theta _{2}}^{\alpha }f\left( \frac{\theta _{1}+\theta _{2}}{2}%
\right) \right] \bigg| \\
&\leq &\frac{4}{p(\alpha p+1)}+\frac{2\left[ \left\vert f^{\prime }(\kappa
_{1})\right\vert ^{q}+\left\vert f^{\prime }(\theta _{2})\right\vert ^{q}%
\right] }{q}
\end{eqnarray*}%
where $p^{-1}+q^{-1}=1$, $\alpha \in \lbrack 0,1]$, $q>1$, $B(\alpha )$ is
normalization function.

\begin{proof}
By using identity that is given in Lemma \ref{ghd65}, we get
\begin{eqnarray*}
&&\bigg|\frac{2(\theta _{2}-\theta _{1})^{\alpha }+(1-\alpha )2^{\alpha
+1}\Gamma (\alpha )}{(\theta _{2}-\theta _{1})^{\alpha +1}}\left[ f(\kappa
_{1})+f(\theta _{2})+2f\left( \frac{\theta _{1}+\theta _{2}}{2}\right) %
\right] \\
&&-\frac{2^{\alpha +1}B(\alpha )\Gamma (\alpha )}{(\theta _{2}-\kappa
_{1})^{\alpha +1}}\left[ ^{AB}I_{\frac{\theta _{1}+\theta _{2}}{2}}^{\alpha
}f(\theta _{1})+_{::::\theta _{1}}^{AB}I^{\alpha }f\left( \frac{\kappa
_{1}+\theta _{2}}{2}\right) \right. \\
&&\left. +_{::\frac{\theta _{1}+\theta _{2}}{2}}^{AB}I^{\alpha }f(\theta
_{2})+^{AB}I_{\theta _{2}}^{\alpha }f\left( \frac{\theta _{1}+\theta _{2}}{2}%
\right) \right] \bigg| \\
&\leq &\int_{0}^{1}(1-\tau _{1})^{\alpha }\left\vert f^{\prime }\left( \frac{%
1+\tau _{1}}{2}\theta _{1}+\frac{1-\tau _{1}}{2}\theta _{2}\right)
\right\vert d\tau _{1}+\int_{0}^{1}\tau _{1}^{\alpha }\left\vert f^{\prime
}\left( \frac{1+\tau _{1}}{2}\theta _{1}+\frac{1-\tau _{1}}{2}\theta
_{2}\right) \right\vert d\tau _{1} \\
&&+\int_{0}^{1}\tau _{1}^{\alpha }\left\vert f^{\prime }\left( \frac{1+\tau
_{1}}{2}\theta _{2}+\frac{1-\tau _{1}}{2}\theta _{1}\right) \right\vert
d\tau _{1}+\int_{0}^{1}(1-\tau _{1})^{\alpha }\left\vert f^{\prime }\left(
\frac{1+\tau _{1}}{2}\theta _{2}+\frac{1-\tau _{1}}{2}\theta _{1}\right)
\right\vert d\tau _{1}.
\end{eqnarray*}%
By using the Young inequality as $xy\leq \frac{1}{p}x^{p}+\frac{1}{q}y^{q}$
\begin{eqnarray*}
&&\bigg|\frac{2(\theta _{2}-\theta _{1})^{\alpha }+(1-\alpha )2^{\alpha
+1}\Gamma (\alpha )}{(\theta _{2}-\theta _{1})^{\alpha +1}}\left[ f(\kappa
_{1})+f(\theta _{2})+2f\left( \frac{\theta _{1}+\theta _{2}}{2}\right) %
\right] \\
&&-\frac{2^{\alpha +1}B(\alpha )\Gamma (\alpha )}{(\theta _{2}-\kappa
_{1})^{\alpha +1}}\left[ ^{AB}I_{\frac{\theta _{1}+\theta _{2}}{2}}^{\alpha
}f(\theta _{1})+_{::::\theta _{1}}^{AB}I^{\alpha }f\left( \frac{\kappa
_{1}+\theta _{2}}{2}\right) \right. \\
&&\left. +_{::\frac{\theta _{1}+\theta _{2}}{2}}^{AB}I^{\alpha }f(\theta
_{2})+^{AB}I_{\theta _{2}}^{\alpha }f\left( \frac{\theta _{1}+\theta _{2}}{2}%
\right) \right] \bigg| \\
&\leq &\frac{1}{p}\int_{0}^{1}(1-\tau _{1})^{\alpha p}d\tau _{1}+\frac{1}{q}%
\int_{0}^{1}\left\vert f^{\prime }\left( \frac{1+\tau _{1}}{2}\theta _{1}+%
\frac{1-\tau _{1}}{2}\theta _{2}\right) \right\vert ^{q}d\tau _{1} \\
&&+\frac{1}{p}\int_{0}^{1}\tau _{1}^{\alpha p}d\tau _{1}+\frac{1}{q}%
\int_{0}^{1}\left\vert f^{\prime }\left( \frac{1+\tau _{1}}{2}\theta _{1}+%
\frac{1-\tau _{1}}{2}\theta _{2}\right) \right\vert ^{q}d\tau _{1} \\
&&+\frac{1}{p}\int_{0}^{1}\tau _{1}^{\alpha p}d\tau _{1}+\frac{1}{q}%
\int_{0}^{1}\left\vert f^{\prime }\left( \frac{1+\tau _{1}}{2}\theta _{2}+%
\frac{1-\tau _{1}}{2}\theta _{1}\right) \right\vert ^{q}d\tau _{1} \\
&&+\frac{1}{p}\int_{0}^{1}(1-\tau _{1})^{\alpha p}d\tau _{1}+\frac{1}{q}%
\int_{0}^{1}\left\vert f^{\prime }\left( \frac{1+\tau _{1}}{2}\theta _{2}+%
\frac{1-\tau _{1}}{2}\theta _{1}\right) \right\vert ^{q}d\tau _{1}.
\end{eqnarray*}%
By using convexity of $|f^{\prime }|^{q}$ and by a simple computation, we
have the desired result.
\end{proof}
\end{thm}

\begin{crly}
In Theorem \ref{9888um}, if we choose $\alpha =1$ we obtain
\begin{eqnarray*}
&&\bigg|\frac{f(\theta _{1})+f(\theta _{2})+2f\left( \frac{\theta
_{1}+\kappa _{2}}{2}\right) }{\theta _{2}-\theta _{1}}-\frac{4}{(\theta
_{2}-\kappa _{1})^{2}}\int_{\theta _{1}}^{\theta _{2}}f(x)dx\bigg| \\
&\leq &\frac{2}{p^{2}+p}+\frac{\left\vert f^{\prime }(\theta
_{1})\right\vert ^{q}+\left\vert f^{\prime }(\theta _{2})\right\vert ^{q}}{q}%
.
\end{eqnarray*}
\end{crly}

\begin{thm}
\label{k?y67f4} Let $f:[\theta _{1},\theta _{2}]\rightarrow \mathbb{R}$ be
differentiable function on $(\theta _{1},\theta _{2})$ with $\kappa
_{1}<\theta _{2}$ and $f^{\prime }\in L_{1}[\theta _{1},\theta _{2}]$. If $%
|f^{\prime }|$ is a concave for $q>1$, then we have
\begin{eqnarray*}
&&\bigg|\frac{2(\theta _{2}-\theta _{1})^{\alpha }+(1-\alpha )2^{\alpha
+1}\Gamma (\alpha )}{(\theta _{2}-\theta _{1})^{\alpha +1}}\left[ f(\kappa
_{1})+f(\theta _{2})+2f\left( \frac{\theta _{1}+\theta _{2}}{2}\right) %
\right] \\
&&-\frac{2^{\alpha +1}B(\alpha )\Gamma (\alpha )}{(\theta _{2}-\kappa
_{1})^{\alpha +1}}\left[ ^{AB}I_{\frac{\theta _{1}+\theta _{2}}{2}}^{\alpha
}f(\theta _{1})+_{::::\theta _{1}}^{AB}I^{\alpha }f\left( \frac{\kappa
_{1}+\theta _{2}}{2}\right) \right. \\
&&\left. +_{::\frac{\theta _{1}+\theta _{2}}{2}}^{AB}I^{\alpha }f(\theta
_{2})+^{AB}I_{\theta _{2}}^{\alpha }f\left( \frac{\theta _{1}+\theta _{2}}{2}%
\right) \right] \bigg| \\
&\leq &\left( \frac{1}{\alpha +1}\right) \Bigg[\left\vert f^{\prime }\left(
\frac{\theta _{1}(\alpha +3)+\theta _{2}(\alpha +1)}{2(\alpha +2)}\right)
\right\vert +\left\vert f^{\prime }\left( \frac{\theta _{1}(2\alpha
+3)+\theta _{2}}{2(\alpha +2)}\right) \right\vert \\
&&+\left\vert f^{\prime }\left( \frac{\theta _{2}(2\alpha +3)+\theta _{1}}{%
2(\alpha +2)}\right) \right\vert +\left\vert f^{\prime }\left( \frac{\kappa
_{2}(\alpha +3)+\theta _{1}(\alpha +1)}{2(\alpha +2)}\right) \right\vert %
\Bigg]
\end{eqnarray*}%
where $\alpha \in \lbrack 0,1]$ and $B(\alpha )$ is normalization function.

\begin{proof}
From Lemma \ref{ghd65} and the Jensen integral inequality, we have
\begin{eqnarray*}
&&\bigg|\frac{2(\theta _{2}-\theta _{1})^{\alpha }+(1-\alpha )2^{\alpha
+1}\Gamma (\alpha )}{(\theta _{2}-\theta _{1})^{\alpha +1}}\left[ f(\kappa
_{1})+f(\theta _{2})+2f\left( \frac{\theta _{1}+\theta _{2}}{2}\right) %
\right]  \\
&&-\frac{2^{\alpha +1}B(\alpha )\Gamma (\alpha )}{(\theta _{2}-\kappa
_{1})^{\alpha +1}}\left[ ^{AB}I_{\frac{\theta _{1}+\theta _{2}}{2}}^{\alpha
}f(\theta _{1})+_{::::\theta _{1}}^{AB}I^{\alpha }f\left( \frac{\kappa
_{1}+\theta _{2}}{2}\right) \right.  \\
&&\left. +_{::\frac{\theta _{1}+\theta _{2}}{2}}^{AB}I^{\alpha }f(\theta
_{2})+^{AB}I_{\theta _{2}}^{\alpha }f\left( \frac{\theta _{1}+\theta _{2}}{2}%
\right) \right] \bigg| \\
&\leq &\left( \int_{0}^{1}(1-\tau _{1})^{\alpha }d\tau _{1}\right)
\left\vert f^{\prime }\left( \frac{\int_{0}^{1}(1-\tau _{1})^{\alpha }\left(
\frac{1+\tau _{1}}{2}\theta _{1}+\frac{1-\tau _{1}}{2}\theta _{2}\right)
d\tau _{1}}{\int_{0}^{1}(1-\tau _{1})^{\alpha }d\tau _{1}}\right)
\right\vert  \\
&&+\left( \int_{0}^{1}\tau _{1}^{\alpha }d\tau _{1}\right) \left\vert
f^{\prime }\left( \frac{\int_{0}^{1}\tau _{1}^{\alpha }\left( \frac{1+\tau
_{1}}{2}\theta _{1}+\frac{1-\tau _{1}}{2}\theta _{2}\right) d\tau _{1}}{%
\int_{0}^{1}\tau _{1}^{\alpha }d\tau _{1}}\right) \right\vert  \\
&&+\left( \int_{0}^{1}\tau _{1}^{\alpha }d\tau _{1}\right) \left\vert
f^{\prime }\left( \frac{\int_{0}^{1}\tau _{1}^{\alpha }\left( \frac{1+\tau
_{1}}{2}\theta _{2}+\frac{1-\tau _{1}}{2}\theta _{1}\right) d\tau _{1}}{%
\int_{0}^{1}\tau _{1}^{\alpha }d\tau _{1}}\right) \right\vert  \\
&&+\left( \int_{0}^{1}(1-\tau _{1})^{\alpha }d\tau _{1}\right) \left\vert
f^{\prime }\left( \frac{\int_{0}^{1}(1-\tau _{1})^{\alpha }\left( \frac{%
1+\tau _{1}}{2}\theta _{2}+\frac{1-\tau _{1}}{2}\theta _{1}\right) d\tau _{1}%
}{\int_{0}^{1}(1-\tau _{1})^{\alpha }d\tau _{1}}\right) \right\vert .
\end{eqnarray*}%
By computing the above integrals we have
\begin{eqnarray*}
&&\bigg|\frac{2(\theta _{2}-\theta _{1})^{\alpha }+(1-\alpha )2^{\alpha
+1}\Gamma (\alpha )}{(\theta _{2}-\theta _{1})^{\alpha +1}}\left[ f(\kappa
_{1})+f(\theta _{2})+2f\left( \frac{\theta _{1}+\theta _{2}}{2}\right) %
\right]  \\
&&-\frac{2^{\alpha +1}B(\alpha )\Gamma (\alpha )}{(\theta _{2}-\kappa
_{1})^{\alpha +1}}\left[ ^{AB}I_{\frac{\theta _{1}+\theta _{2}}{2}}^{\alpha
}f(\theta _{1})+_{::::\theta _{1}}^{AB}I^{\alpha }f\left( \frac{\kappa
_{1}+\theta _{2}}{2}\right) \right.  \\
&&\left. +_{::\frac{\theta _{1}+\theta _{2}}{2}}^{AB}I^{\alpha }f(\theta
_{2})+^{AB}I_{\theta _{2}}^{\alpha }f\left( \frac{\theta _{1}+\theta _{2}}{2}%
\right) \right] \bigg| \\
&\leq &\left( \frac{1}{\alpha +1}\right) \Bigg[\left\vert f^{\prime }\left(
\frac{\theta _{1}(\alpha +3)+\theta _{2}(\alpha +1)}{2(\alpha +2)}\right)
\right\vert +\left\vert f^{\prime }\left( \frac{\theta _{1}(2\alpha
+3)+\theta _{2}}{2(\alpha +2)}\right) \right\vert  \\
&&+\left\vert f^{\prime }\left( \frac{\theta _{2}(2\alpha +3)+\theta _{1}}{%
2(\alpha +2)}\right) \right\vert +\left\vert f^{\prime }\left( \frac{\kappa
_{2}(\alpha +3)+\theta _{1}(\alpha +1)}{2(\alpha +2)}\right) \right\vert %
\Bigg].
\end{eqnarray*}%
So, the proof is completed.
\end{proof}
\end{thm}

\begin{crly}
In Theorem \ref{k?y67f4}, if we choose $\alpha =1$ we obtain
\begin{eqnarray*}
&&\bigg|\frac{f(\theta _{1})+f(\theta _{2})+2f\left( \frac{\kappa
_{1}+\theta _{2}}{2}\right) }{\theta _{2}-\theta _{1}}-\frac{4}{(\kappa
_{2}-\theta _{1})^{2}}\int_{\theta _{1}}^{\theta _{2}}f(x)dx\bigg| \\
&\leq &\left( \frac{1}{4}\right) \Bigg[\left\vert f^{\prime }\left( \frac{%
2\theta _{1}+\theta _{2}}{3}\right) \right\vert +\left\vert f^{\prime
}\left( \frac{5\theta _{1}+\theta _{2}}{6}\right) \right\vert \\
&&+\left\vert f^{\prime }\left( \frac{5\theta _{2}+\theta _{1}}{6}\right)
\right\vert +\left\vert f^{\prime }\left( \frac{2\theta _{2}+\theta _{1}}{3}%
\right) \right\vert \Bigg].
\end{eqnarray*}
\end{crly}

\begin{thm}
\label{ght56d} Let $f:[\theta _{1},\theta _{2}]\rightarrow \mathbb{R}$ be
differentiable function on $(\theta _{1},\theta _{2})$ with $\kappa
_{1}<\theta _{2}$ and $f^{\prime }\in L_{1}[\theta _{1},\theta _{2}]$. If $%
|f^{\prime }|^{q}$ is a concave function, we have
\begin{eqnarray*}
&&\bigg|\frac{2(\theta _{2}-\theta _{1})^{\alpha }+(1-\alpha )2^{\alpha
+1}\Gamma (\alpha )}{(\theta _{2}-\theta _{1})^{\alpha +1}}\left[ f(\kappa
_{1})+f(\theta _{2})+2f\left( \frac{\theta _{1}+\theta _{2}}{2}\right) %
\right] \\
&&-\frac{2^{\alpha +1}B(\alpha )\Gamma (\alpha )}{(\theta _{2}-\kappa
_{1})^{\alpha +1}}\left[ ^{AB}I_{\frac{\theta _{1}+\theta _{2}}{2}}^{\alpha
}f(\theta _{1})+_{::::\theta _{1}}^{AB}I^{\alpha }f\left( \frac{\kappa
_{1}+\theta _{2}}{2}\right) \right. \\
&&\left. +_{::\frac{\theta _{1}+\theta _{2}}{2}}^{AB}I^{\alpha }f(\theta
_{2})+^{AB}I_{\theta _{2}}^{\alpha }f\left( \frac{\theta _{1}+\theta _{2}}{2}%
\right) \right] \bigg| \\
&\leq &\frac{2}{(\alpha p+1)^{\frac{1}{p}}}\left[ \left\vert f^{\prime
}\left( \frac{3\theta _{1}+\theta _{2}}{4}\right) \right\vert +\left\vert
f^{\prime }\left( \frac{3\theta _{2}+\theta _{1}}{4}\right) \right\vert %
\right]
\end{eqnarray*}%
where $p^{-1}+q^{-1}=1$, $\alpha \in \lbrack 0,1]$, $q>1$.

\begin{proof}
By using the Lemma \ref{ghd65} and H\"{o}lder integral inequality, we can
write
\begin{eqnarray*}
&&\bigg|\frac{2(\theta _{2}-\theta _{1})^{\alpha }+(1-\alpha )2^{\alpha
+1}\Gamma (\alpha )}{(\theta _{2}-\theta _{1})^{\alpha +1}}\left[ f(\kappa
_{1})+f(\theta _{2})+2f\left( \frac{\theta _{1}+\theta _{2}}{2}\right) %
\right] \\
&&-\frac{2^{\alpha +1}B(\alpha )\Gamma (\alpha )}{(\theta _{2}-\kappa
_{1})^{\alpha +1}}\left[ ^{AB}I_{\frac{\theta _{1}+\theta _{2}}{2}}^{\alpha
}f(\theta _{1})+_{::::\theta _{1}}^{AB}I^{\alpha }f\left( \frac{\kappa
_{1}+\theta _{2}}{2}\right) \right. \\
&&\left. +_{::\frac{\theta _{1}+\theta _{2}}{2}}^{AB}I^{\alpha }f(\theta
_{2})+^{AB}I_{\theta _{2}}^{\alpha }f\left( \frac{\theta _{1}+\theta _{2}}{2}%
\right) \right] \bigg| \\
&\leq &\left( \int_{0}^{1}(1-\tau _{1})^{\alpha p}d\tau _{1}\right) ^{\frac{1%
}{p}}\left( \int_{0}^{1}\left\vert f^{\prime }\left( \frac{1+\tau _{1}}{2}%
\kappa _{1}+\frac{1-\tau _{1}}{2}\theta _{2}\right) \right\vert ^{q}d\tau
_{1}\right) ^{\frac{1}{q}} \\
&&+\left( \int_{0}^{1}\tau _{1}^{\alpha p}d\tau _{1}\right) ^{\frac{1}{p}%
}\left( \int_{0}^{1}\left\vert f^{\prime }\left( \frac{1+\tau _{1}}{2}\theta
_{1}+\frac{1-\tau _{1}}{2}\theta _{2}\right) \right\vert ^{q}d\tau
_{1}\right) ^{\frac{1}{q}} \\
&&+\left( \int_{0}^{1}\tau _{1}^{\alpha p}d\tau _{1}\right) ^{\frac{1}{p}%
}\left( \int_{0}^{1}\left\vert f^{\prime }\left( \frac{1+\tau _{1}}{2}\theta
_{2}+\frac{1-\tau _{1}}{2}\theta _{1}\right) \right\vert ^{q}d\tau
_{1}\right) ^{\frac{1}{q}} \\
&&+\left( \int_{0}^{1}(1-\tau _{1})^{\alpha p}d\tau _{1}\right) ^{\frac{1}{p}%
}\left( \int_{0}^{1}\left\vert f^{\prime }\left( \frac{1+\tau _{1}}{2}\theta
_{2}+\frac{1-\tau _{1}}{2}\theta _{1}\right) \right\vert ^{q}d\tau
_{1}\right) ^{\frac{1}{q}}.
\end{eqnarray*}

By using concavity of $|f^{\prime }|^{q}$ and Jensen integral inequality, we
get
\begin{eqnarray*}
&&\int_{0}^{1}\left\vert f^{\prime }\left( \frac{1+\tau _{1}}{2}\theta _{1}+%
\frac{1-\tau _{1}}{2}\theta _{2}\right) \right\vert ^{q}d\tau _{1} \\
&=&\int_{0}^{1}\tau _{1}^{0}\left\vert f^{\prime }\left( \frac{1+\tau _{1}}{2%
}\theta _{1}+\frac{1-\tau _{1}}{2}\theta _{2}\right) \right\vert ^{q}d\tau
_{1} \\
&\leq &\left( \int_{0}^{1}\tau _{1}^{0}d\tau _{1}\right) \left\vert
f^{\prime }\left( \frac{\int_{0}^{1}\tau _{1}^{0}\left( \frac{1+\tau _{1}}{2}%
\theta _{1}+\frac{1-\tau _{1}}{2}\theta _{2}\right) d\tau _{1}}{%
\int_{0}^{1}\tau _{1}^{0}d\tau _{1}}\right) \right\vert ^{q} \\
&=&\left\vert f^{\prime }\left( \frac{3\theta _{1}+\theta _{2}}{4}\right)
\right\vert ^{q}.
\end{eqnarray*}%
Similarly
\begin{equation*}
\int_{0}^{1}\left\vert f^{\prime }\left( \frac{1+\tau _{1}}{2}\theta _{2}+%
\frac{1-\tau _{1}}{2}\theta _{1}\right) \right\vert ^{q}d\tau _{1}\leq
\left\vert f^{\prime }\left( \frac{3\theta _{2}+\theta _{1}}{4}\right)
\right\vert ^{q}
\end{equation*}%
so, we obtain
\begin{eqnarray*}
&&\bigg|\frac{2(\theta _{2}-\theta _{1})^{\alpha }+(1-\alpha )2^{\alpha
+1}\Gamma (\alpha )}{(\theta _{2}-\theta _{1})^{\alpha +1}}\left[ f(\kappa
_{1})+f(\theta _{2})+2f\left( \frac{\theta _{1}+\theta _{2}}{2}\right) %
\right] \\
&&-\frac{2^{\alpha +1}B(\alpha )\Gamma (\alpha )}{(\theta _{2}-\kappa
_{1})^{\alpha +1}}\left[ ^{AB}I_{\frac{\theta _{1}+\theta _{2}}{2}}^{\alpha
}f(\theta _{1})+_{::::\theta _{1}}^{AB}I^{\alpha }f\left( \frac{\kappa
_{1}+\theta _{2}}{2}\right) \right. \\
&&\left. +_{::\frac{\theta _{1}+\theta _{2}}{2}}^{AB}I^{\alpha }f(\theta
_{2})+^{AB}I_{\theta _{2}}^{\alpha }f\left( \frac{\theta _{1}+\theta _{2}}{2}%
\right) \right] \bigg| \\
&\leq &\frac{2}{(\alpha p+1)^{\frac{1}{p}}}\left[ \left\vert f^{\prime
}\left( \frac{3\theta _{1}+\theta _{2}}{4}\right) \right\vert +\left\vert
f^{\prime }\left( \frac{3\theta _{2}+\theta _{1}}{4}\right) \right\vert %
\right] .
\end{eqnarray*}
\end{proof}
\end{thm}

\begin{crly}
In Theorem \ref{ght56d}, if we choose $\alpha =1$ we obtain
\begin{eqnarray*}
&&\bigg|\frac{f(\theta _{1})+f(\theta _{2})+2f\left( \frac{\kappa
_{1}+\theta _{2}}{2}\right) }{\theta _{2}-\theta _{1}}-\frac{4}{(\kappa
_{2}-\theta _{1})^{2}}\int_{\theta _{1}}^{\theta _{2}}f(x)dx\bigg| \\
&\leq &\frac{1}{(p+1)^{\frac{1}{p}}}\Bigg[\left\vert f^{\prime }\left( \frac{%
3\theta _{1}+\theta _{2}}{4}\right) \right\vert +\left\vert f^{\prime
}\left( \frac{3\theta _{2}+\theta _{1}}{4}\right) \right\vert \Bigg].
\end{eqnarray*}
\end{crly}

\section{Conclusion}

In this study, an integral identity including Atangana-Baleanu integral operators has been proved. Some integral inequalities are established by using H\"older inequality, Power-mean inequality, Young inequality and convex functions with the help of Lemma 2.1  which has the potential to produce Bullen type inequalities. Some special cases of the results in this general form have been pointed out. Researchers can establish new equations such as the integral identity in the study and reach similar inequalities of these equality-based inequalities.

\section{Acknowledgments}

The publication has been prepared with the support of GNAMPA 2019 and the RUDN University Strategic Academic Leadership Program.

\section*{Funding}

GNAMPA 2019 and the RUDN University Strategic Academic Leadership Program.

\section*{Availability of data and materials}

Data sharing is not applicable to this paper as no datasets were generated
or analyzed during the current study.

\section*{Competing interests}

The authors declares that there is no conflict of interests regarding the
publication of this paper.

\section*{Author's contributions}

All authors jointly worked on the results and they read and approved the
final manuscript.

\end{document}